\documentstyle[12pt,a4,amscd,amssymb,verbatim]{amsart}

\newtheorem{Thm}{Theorem}

\newtheorem{Corollary}{Corollary}  
\newtheorem{Lemma}{Lemma}

\newtheorem{Definition}[Thm]{Definition}

\newcommand{\Ker}{\operatorname{Ker}}

\newcommand{\Hom}{\operatorname{Hom}}
\newcommand{\End}{\operatorname{End}}
\newcommand{\Ext}{\operatorname{Ext}}
\newcommand{\Extder}{\operatorname{{ \mathit{Ext} }}}
\newcommand{\modc}{\operatorname{mod}}
\newcommand{\Modc}{\operatorname{Mod}}
\newcommand{\mof}{\operatorname{mof}}
\newcommand{\Mof}{\operatorname{Mof}}
\newcommand{\Homder}{\operatorname{{ \mathit{Hom} }}}
\newcommand{\Endder}{\operatorname{{ \mathit{End} }}}
\newcommand{\smallend}{\operatorname{{ \mathit{end} }}}
\newcommand{\RHom}{\operatorname{{ \mathit{RHom} }}}

\begin{document}

\title{Koszul duality of translation---and Zuckerman functors} 
\author{Steen Ryom-Hansen}
\address{ Matematisk Afdeling Universitetsparken 5 DK-2100 K{\o}benhavn 
{\O} Danmark steen@@math.ku.dk }
\footnote{Supported by EPSRC grant M22536 and by the TMR-network
algebraic Lie Theory ERB-FMRX-CT97-0100.} 

\maketitle

\begin{abstract}
We review Koszul duality in representation theory of 
category $ \cal O $, especially 
we give a new presentation of the Koszul duality functor.
Combining this with work of Backelin, we show that the translation 
and Zuckerman functors are Koszul dual to each other, thus verifying a 
conjecture of Bernstein, Frenkel and Khovanov. Finally we use Koszul 
duality to give a short proof of the Enright-Shelton equivalence.
\end{abstract}

\section {\bf Introduction. }

\noindent
The definite exposition of Koszul duality in representation theory is the 
paper of Beilinson, Ginzburg and Soergel [5]. The main theme 
of that paper is that the 
category $ \cal O $ of Bernstein, Gelfand and Gelfand is the module category
of a Koszul ring, from which many properties of $ \cal O $ can be deduced very
elegantly. They also consider singular as well as parabolic versions of
$ \cal O $ and
show that they are are equivalent to each other under the Koszul duality 
functor. Two of the main tools to obtain these results were the Zuckerman 
functors and the translation functors.

\medskip
\noindent

In a recent paper by Bernstein, Frenkel and Khovanov [6] these functors are 
studied on certain simultaneously singular and parabolic subcategories of 
$ \cal O $, in type A. They show, partially conjecture, that each 
one of them 
-- in combination with an equivalence of certain categories due to Enright-Shelton \nolinebreak --
can be 
used to categorify the Temperley-Lieb algebra and also conjecture that these
two pictures should be Koszul dual to each other.

\medskip
\noindent

The purpose of this note is to show that the translation-- and Zuckerman 
functors are indeed Koszul dual to each other. By this we mean that both admit 
graded versions and that these correspond under the Koszul duality 
functor. Although this of course is one of the main philosophical points of 
[5], a rigorous proof was never given.

\medskip
\noindent

We furthermore use the Koszul duality theory to 
give a simple proof of the Enright-Shelton equivalence. 

\medskip

It is a pleasure to thank Wolfgang Soergel for many useful
discussions. Most of this 
work was done while I was a research assistant at the University 
of Freiburg, Germany. 

\medskip

Added: Recently C. Stroppel [11] has obtained a proof of the full
Bernstein-Frenkel-Khovanov conjectures, in part using results from the present paper.

\section{\bf Preliminaries}

\noindent
In this section we will mostly recall some of the basic definitions 
and concepts 
from the theory of Koszul duality, see also [5].

\medskip

Let $ \mathfrak g $ be a complex semisimple Lie algebra containing a Borel
algebra $ \mathfrak b $ with Cartan part $ \mathfrak h $ and let $ \cal O $ 
be the category of $ \mathfrak g $-modules associated to it by Bernstein, 
Gelfand and Gelfand. For $ \lambda \, \in \mathfrak h ^* $ integral but 
possibly singular, denote by $ \cal O_{\lambda} $ the subcategory of 
$ \cal O $ 
consisting of modules of generalized central character $ \chi_{\lambda} $ --
we refer to it as the singular category $ \cal O $. Let $ S $ be the set of 
simple reflections corresponding to our data, and let $ S_{\lambda } $ be the 
subset consisting of those reflections that fix $ \lambda $ under the dot 
operation. Then $ S_{\lambda } $ defines a parabolic subalgebra
$ \mathfrak q ( \lambda )$. We define the parabolic category 
$ \cal O^{\lambda} $ to consist of the $ \mathfrak q ( \lambda )$-finite 
objects in $ \cal O $. For $ \lambda = 0 $ we omit the index, i.e. we write 
$ {\cal O}_0 \, = \, \cal O $, this should not cause confusion. 

\medskip
\noindent

Let $ P $ denote the sum of all indecomposable projectives in $ \cal O $, hence
$ P $ is a projective generator of $ \cal O $. Analogously, we can construct 
projective generators $ P_{\lambda} $ (resp. $ P^{\lambda} $) of 
$ \cal{ O}_{\lambda} $ (resp. $ {\cal O}^{\lambda} $). 
Let $ A  =  \End_{\cal O} P $ (resp. $  A_{\lambda} = 
\End_{{\cal O }_{\lambda}} P_{\lambda}, \, \, 
 A^{\lambda}  =  \End_{{\cal O }^{\lambda}} P^{\lambda} $). 
By general principles, we can then identify $ \cal O $ with 
$ { \Modc }{\textstyle-}A $,
the category of finitely generated $ A $-right modules (and analogously for 
$ A_{\lambda} $, $ A^{\lambda} $).

\medskip
\noindent

Let $ T_0^{\lambda}:\, \cal{O}\, \rightarrow \cal{O}_{\lambda}\,\, $,
$ T^0_{\lambda}:\, \cal{O}_{\lambda} \rightarrow \cal{O}\,   $
be the Jantzen translation functors onto and out of the wall. 
Passing to $ {\Modc} { \textstyle -}A $, 
$ T_0^{\lambda} $ corresponds to 
the functor 
$$ { \Modc}{\textstyle -} A \, \rightarrow 
{ \Modc } {\textstyle - } A_{\lambda}:\, \, M \, \mapsto 
M \otimes_{A}  X $$
where $ X $ is the $ ( A,A_{\lambda} ) $ bimodule 
$ X \, = \, \Hom_{\cal O}( \, P_{\lambda}, T_0^{\lambda} \, P \,) $.
There is a similar description of 
$ T^0_{\lambda}:\, \cal{O}_{\lambda} \rightarrow \cal{O} $.

\medskip
\noindent

Let $ \tau_{\lambda}: \cal O \, \rightarrow \, \cal {O}^{\lambda} $ be 
the parabolic truncation functor, by definition it 
takes $ M \, \in \cal O $ 
to its largest $ \mathfrak {q}( \lambda ) $-finite quotient. It is 
right exact and left  
adjoint to the inclusion functor $ \iota_{\lambda} : {\cal O }^{\lambda} \, 
\rightarrow \, \cal O $, which is exact, and it thus  takes projectives to 
projectives. It even takes indecomposable projectives to indecomposable 
projectives.
Its top degree left derived functor is the Zuckerman functor 
that takes a module $ M \,\in {\cal O } $ to its largest 
$ \mathfrak {q}( \lambda ) $-finite submodule. 

\medskip
\noindent

On the $ { \Modc}{\textstyle -}A $ level, $ \tau_{\lambda} $ 
is given by the tensor product 
with the $ ( A,A^{\lambda} ) $-bimodule 
$ Y  = \Hom_{{\cal O}^{\lambda}}
( \, P^{\lambda},\, \tau_{\lambda} P\, ) $. However, from the above 
considerations we deduce that $ Y = A^{\lambda} $, the left $ A $-structure 
coming from $ \tau_{\lambda} $ and the right $ A^{\lambda} $-structure 
coming from the 
multiplication in $ A^{\lambda} $.

\medskip
\noindent
Let us finally recall the definition of a Koszul ring. 
\medskip

\begin{Definition} Let $R = \bigoplus_{i \geq 0}R_i$ be a 
positively graded ring with $R_0$ semisimple and put
$R_+ = \bigoplus_{i > 0} R_i$. $R$ is called Koszul 
if the right module $R_0 \cong R/R_+$\
admits a graded projective resolution $P^{\bullet} 
\twoheadrightarrow R_0$ such
that $P^{i}$ 
is generated by its component $P^{i}_i$ in 
degree $i$. 
The Koszul dual ring of $ R $ is defined as 
$ R^!:= \Ext^{\bullet}_R( R_0,R_0) $.
\end{Definition} 

\medskip

\noindent
The main point of [5] is now that all the rings appearing above can be 
given a Koszul grading. 

\section{\bf Koszul duality. }

\noindent
The main purpose of this section is to give a new construction of the 
Koszul duality functor. 

\medskip
Koszul duality is an equivalence at the level of derived categories between
the graded module categories of a Koszul ring and its dual. 
It exists under some 
mild finiteness conditions. The concept appeared for the first time in the 
paper [3] 
where it is shown that for any finite dimensional vector space $ V $ the 
derived graded module categories of the 
symmetric algebra $ SV $ and the exterior algebra $ \bigwedge V^* $ are 
equivalent. The argument used there carries over to general Koszul rings and 
is the one used in [5].

\medskip 
\noindent

We shall here present another approach to the Koszul duality functor, 
using the language of differential graded algebras (DG-algebras).
Although also the papers [4] and [9] link DG-algebras with 
Koszul duality, 
we were rather inspired by the construction of the localization functor in 
the book of Bernstein and Lunts [7].

\medskip 

Let us start out by repeating the basic definitions.

\begin{Definition}
A DG-algebra $ {\cal A} = ( A,d ) $ is a graded associative algebra 
$ A = \bigoplus_{i = - \infty }^{\infty} A^i $ with a unit $ 1_A \in A^0 $
and an 
additive endomorphism $ d $ of degree 1, s.t.
$$ d^2 = 0 $$ $$ d( a \cdot b )= da \cdot b + ( -1)^{deg(a)} a \cdot db $$
$$  and \,\,\, d( 1_A ) = 0. $$
\end{Definition}

\begin{Definition} 
A left module ${\cal M}= ( M,d_M ) $ over a DG-algebra  
$ {\cal A} = ( A,d ) $ is a 
graded unitary right $ A $-module 
$  M = \bigoplus_{i = - \infty }^{\infty} M^i $ with an additive
endomorphism
$ d_M:{\cal M}  \rightarrow {\cal M} $ of degree 1, s.t. $ d_M^2 = 0 $ and 
$$ d_M(am) = da \cdot m + ( -1)^{deg(a)} a\cdot d_Mm $$
\end{Definition}

\begin{Definition} 
A right module ${\cal M}= ( M,d_M ) $ over a DG-algebra  
$ {\cal A} = ( A,d ) $ 
is a graded unitary right $ A $-module 
$  M  = \bigoplus_{i = - \infty }^{\infty} M^i $ with an additive 
endomorphism
$ d_M: {\cal M} \rightarrow {\cal M} $ of degree 1, s.t. $ d_M^2 = 0 $ and 
$$ d_M(ma) = d_M m \cdot a + ( -1)^{deg(m)}m \cdot da $$
\end{Definition}

\medskip

{ \it Unless otherwise stated, $ A $ will from now on be the ring $ \End _{\cal O } (P)$ from the former section. }

\medskip

By the selfduality theorem of Soergel [10]
we know that $ A = \Ext^{\bullet}_{\cal O } ( L, L ) $, where $ L $ is the 
sum of all simples in $ \cal{O} $. This provides $ A $ with a grading 
which by [10] or [5] is a Koszul grading.

\medskip

We regard $ A $ as a DG-algebra $ \cal{A}$ with $ A = A^0 $. At this stage 
we neglect the grading on $ A $, hence $ \cal{A}$ just has one grading, the 
DG-algebra grading.

\medskip

Let $ K^{\bullet} $ be the Koszul
complex of $ A $, see e.g. [5]. It is a projective resolution of 
$ k $ (which corresponds to the sum of all simples in
$ \cal O $). In [5] $ K^{\bullet} $ is a resolution of 
left $ A $-modules of 
$ k $; we however prefer to modify it to a resolution of 
right $ A $-modules. Furthermore we change the sign of the indices so 
that the differential has degree $ + 1 $.

\medskip
Let $ k_w $ denote the simple $ A $-module corresponding to 
$ L_w \in \cal O $ and  
let $ K^{\bullet}_{w} $ be a projective resolution 
of $ k_{w} $.
The Koszul complex 
$K^{\bullet}$ is then quasiisomorphic to $ \bigoplus_{\lambda}
K^{\bullet}_w $.

\medskip

We can consider $ K^{\bullet} $ as a DG right module of 
the DG-algebra 
$\cal A$, and as such it is $ \cal{K} $-projective in the sense of [7]: the 
homotopy category of DG-modules of $ A $ is simply the homotopy category of 
complexes of $ A $-modules, hence any quasiisomorphism from $ K^{\bullet}  $
to an DG-module $ \cal M $ is a homotopy isomorphism. 

\medskip

Given two DG-modules $\cal M, N $ of the (arbitrary) DG-algebra $ \cal A $,
recall the construction of the complex 
$ \Homder_{\cal A}^{\bullet}( M,N ) $: 
$$ \Homder^n_{\cal A }( M,N ) := \Hom_A ( M, N[n] ) $$
$$ df:= d_Nf - (-1)^n f d_M, \, \, f \in \Homder_A^n( M,N ). $$

If we perform this construction on $K^{\bullet} $, we obtain a DG-algebra 
$$ \Endder_{\cal A }^{\bullet}( K^{\bullet} ) = 
\Homder^{\bullet}_{\cal A}( K^{\bullet} ,K^{\bullet} ),$$ the multiplication 
being given by composition of maps
(notice that this construction does not involve the Koszul 
grading on $ A $).

The signs match up to give
$K^{\bullet}$ the structure of an $ \Endder_{\cal A }^{\bullet}( K^{\bullet} )$ 
left DG-module.

\medskip
\noindent

Let for any DG-algebra $ \cal A $,  
$ {\cal D}^b {\Modc -}\cal A ) $ denote the bounded derived category of 
right $ \cal A $-modules, see for example Bernstein-Lunts [7] for an 
exposition of this theory.

\medskip
\noindent

Let $ {\cal D}^b ( {\Mof-}  A  ) $ denote the bounded derived category 
of finitely generated $ A $-modules, and let 
$ {\cal D}^b ( {\Mof-} {\cal  A}  ) $ be the corresponding category of 
$ \cal A $-modules. 
Let $ {\cal D}^b( {\Mof-}  \Endder_{\cal A }^{\bullet}( K^{\bullet} )) $ be
the subcategory of 
$ {\cal D}^b( {\Modc-}  \Endder_{\cal A }^{\bullet}( K^{\bullet} )) $ 
consisting of finitely generated $\Endder_{\cal A }^{\bullet}( K^{\bullet} ) $-modules in the algebra sense.
Then we have

\medskip

\begin{Lemma}
${\cal D}^b( {\Mof-}  \Endder_{\cal A }^{\bullet}( K^{\bullet} )) $ is
generated, in the sense of triangulated categories, by the $\Homder_{\cal A
}^{\bullet}( K^{\bullet},K^{\bullet}_w) $'s.
\end{Lemma}
\begin{pf*}{Proof} 
Since our ring $ A \, $ is Koszul, the DG-algebras 
$ \Endder_{\cal A }^{\bullet}\left( K^{\bullet}\right)$ and
$ ( \Extder_{ A }^{\bullet}(k,k ), d=0 )$ are 
quasiisomorphic [5, Theorem 2.10.1] and hence their bounded derived categories 
are equivalent. Under this equivalence the
$ \Endder_{\cal A }^{\bullet}\left( K^{\bullet}\right)$-module 
$\Homder_{\cal A}^{\bullet}( K^{\bullet},K^{\bullet}_w) $
becomes the 
$ (\Extder_{\cal A }^{\bullet}\left( k,k  \right), d=0)$-module 
$ (\Extder_{\cal A }^{\bullet}\left( k, k_w \right), d=0)$.

\medskip 

Now a certain polynomial DG-algebra $ \cal B $ with $d=0$ is studied in chapter 11 of [7]. 
Corollary 11.1.5 there is the statement that 
${\cal D}^b( {\Mof-} {\cal B}) $ is
generated by $ \cal B $ itself. The proof of this 
Corollary 11.1.5 relies on the differential being zero for $ \cal B $ and of
the fact that finitely generated projective modules over a
polynomial algebra are free. We do not have this last property in 
our situation and should therefore adjust the statement of the Lemma accordingly.
Let us give a few details on the modifications: 

\medskip 



Let thus $ ({\cal M}, d_M) $  be an object of
${\cal D}^b( {\Mof-}  (\Extder_{\cal A }^{\bullet}( k,k ), d=0)) $.
Then we have a triangle in that category of the form 
$$  \Ker d_M 
\rightarrow  M \rightarrow M/\Ker d_M  \rightarrow $$
But the first and the third terms have zero differential so we may 
assume that $ d_M = 0 $. 

\medskip 

Let now 
$$ 0 \rightarrow  P^{-n} \rightarrow P^{-n+1} \rightarrow
\ldots \stackrel{\epsilon} {\rightarrow} M \rightarrow 0$$
be a projective graded resolution of $M $ considered as 
$\Extder_{\cal A }^{\bullet}( k,k ) $-module.
Defining $ {\cal P} \in  
{\cal D}^b( {\Mof-}  (\Extder_{\cal A }^{\bullet}( k,k ), d=0))$
as having 
$k$'th DG-part $ \oplus_i P^i_{k-i} $ and differential from 
the resolution, we get a quasiisomorphism 
$  \cal P \cong M $ induced by  $ \epsilon $.

\medskip 
On the other hand the finitely generated, graded projectives over 
$ \Extder_{\cal A }^{\bullet}\left( k,k  \right)$
are the direct sums of the
$ \Extder_{\cal A }^{\bullet}\left(  k, k_w   \right)$'s 
and we are done.

\end{pf*} 

\medskip

Now since $ K^{\bullet} $ is $ \cal K $-projective, it can be used in itself 
to calculate \hfill\break
$ \RHom_{\cal A }(K^{\bullet}, M ) $ for an $ \cal A $-module
$\cal M $. But $ K^{\bullet} $ is an 
$ \Endder_{\cal A }^{\bullet}( K^{\bullet} )$-
module, so we get a functor:
$$  \RHom_{\cal A }(K^{\bullet}, - ):  {\cal D}^b ({\Modc-}{\cal A} ) 
\rightarrow {\cal D}^b ( {\Modc-}  \Endder_{\cal A }^{\bullet}( K^{\bullet} ))$$

This might require a little consideration: one should check that \hfill\break
$ \RHom_{\cal A }(K^{\bullet}, - ) $ takes homotopies to homotopies, even when 
considered as a functor to the category of right 
$ \Endder_{\cal A }^{\bullet}( K^{\bullet}) $-modules, and likewise for 
quasiisomorphisms. For homotopies, one checks this by hand (after having 
consulted e.g. Bernstein-Lunts [7] for the definition of the 
homotopy category of 
DG-modules). For quasiisomorphisms one uses the standard argument: if 
$ f: M \rightarrow N $ is a quasiisomorphism, then 
$ \Homder_{\cal A }^{\bullet}( K^{\bullet}, C(f) )$ is acyclic ($ C(f) $ 
denoting  the cone of $ f $) and so on: the action of  
$ \Endder_{\cal A }^{\bullet}( K^{\bullet} ) $ is irrelevant for this argument.

\medskip
\noindent

Now we also have a functor in the other direction:
$$ -\stackrel{L}{\otimes}_{ \Endder_{\cal A }^{\bullet}( K^{\bullet} ) }
K^{\bullet}
: \, \,
{\cal D}^b ({\Modc-}{\cal A} ) 
\leftarrow {\cal D}^b ( {\Modc-}  \Endder_{\cal A }^{\bullet}( K^{\bullet} )) $$

Recall that this construction involves the replacement of $${\cal N} \in 
{\cal D}^b ( {\Modc-}  \Endder_{\cal A }^{\bullet}( K^{\bullet} ))$$ 
by a $ \cal K $-flat object (actually a $ \cal K $-projective object will 
do), to which it is quasiisomorphic.
The action of $ A $ on $ K^{\bullet} $ then provides 
$ N \stackrel{L}{\otimes}_{ \Endder_{\cal A }^{\bullet}( K^{\bullet} )}
K^{\bullet}   $ with the structure of DG-module 
over $ \cal A $: 
$$ \begin{array}{rcl}
d(n \otimes ka) & = & d(n) \otimes ka + (-1)^{deg(n)} n \otimes d(ka) \\
                & = & d(n) \otimes ka + (-1)^{deg(n)} n \otimes d(k)a  \\
                & = & d(n \otimes k)a 
\end{array} $$ 
since the differential of $ \cal A $ is zero. Once again, one should  
here check that the functor makes sense as a functor into the category of 
$ \cal A $-modules. 

\medskip
\noindent

We can then prove the following Theorem:

\medskip

\begin{Thm}
The functor 
$ \RHom_{\cal A }(K^{\bullet}, - )$  
establishes an equivalence of the triangulated categories 
$ {\cal D}^b (\Mof- {\cal A} )   $ and 
$ {\cal D}^b ( {\Mof-}  \Endder_{\cal A }^{\bullet}\left( K^{\bullet}) 
\right)$. 
The inverse functor is  
$ -\stackrel{L}{\otimes}_{ \Endder_{\cal A }^{\bullet}( K^{\bullet} ) }
K^{\bullet}.$ 
\end{Thm}

\begin{pf*}{Proof} 
Note first that 
$ {\cal D}^b ( {\Mof-}  A  ) $ is generated by 
the $ K_w^{\bullet}  $, 
since they are resolutions of the simple $A$-modules 
and short exact sequences $\cal O$-modules induce 
triangles in $ {\cal D}^b ( {\Modc-}  A)  $. 

\medskip
But then 
$ \RHom_{\cal A }(K^{\bullet}, M ) $ belongs to
${\cal D}^b( {\Mof-}  \Endder_{\cal A }^{\bullet}( K^{\bullet} )) $
for 
$ M \in {\cal D}^b( {\Mof-}  \Endder_{\cal A }^{\bullet}( K^{\bullet} )) $, 
since it clearly does so for each of the generators  $  K_w^{\bullet} $.

\medskip


We furthermore see that 
$ \RHom_{\cal A }(K^{\bullet}, M ) $ is a $ \cal K $-projective 
$  \Endder_{\cal A }^{\bullet}( K^{\bullet}) $-module for any $ M $ in 
$ {\cal D}^b ( {\Mof-} {\cal  A}  ) $
since it 
clearly holds for $  K_w^{\bullet} $ and the $\cal K $-projectives form a 
triangulated subcategory of  
$ {\cal K}( {\Modc-}  \Endder_{\cal A }^{\bullet}( K^{\bullet} )) $, the 
homotopy category of $  \Endder_{\cal A }^{\bullet}( K^{\bullet} ) $-modules.

\medskip
\noindent

We then obtain a  
morphism $ \psi_M $ in $ {\cal D}^b( {\Modc}-\cal A ) $ in the 
following way.
$$ \psi_M : \, \RHom_{\cal A }(K^{\bullet}, M )
\stackrel{L}{\otimes}_{{ \Endder_{\cal A }^{\bullet}}( K^{\bullet} )} K^{\bullet} \rightarrow 
M: \, \, \, f \otimes k \mapsto f(k) $$

One checks that $ \psi_M $ defines a natural transformation from the functor 
$$ \RHom_{\cal A }(K^{\bullet}, - )
\stackrel{L}{\otimes}_{{ \Endder_{\cal A }^{\bullet}}( K^{\bullet} )}
K^{\bullet} $$ to the identity functor $ Id $. 
It is a quasiisomorphism for $ M = K_w^{\bullet} $, and thus for all $ M $.

\medskip

On the other hand we have for a $ \cal K $-projective 
$ \Endder_{\cal A }^{\bullet}( K^{\bullet} ) $-module $ N $ the canonical morphism 
$ \phi_N $ given as follows:

$$ \begin{array}{rcl} \phi_N : \, \,
N & \rightarrow & \RHom_{\cal A }\left( K^{\bullet}, N \stackrel{L}
{\otimes}_{ \Endder_{\cal A }^{\bullet}( K^{\bullet})}  K^{\bullet} \right ) \\
 n_i \in N^i  & \mapsto & (\,  k_j \in K^j \mapsto n_i \otimes k_j  \,  ) 
\end{array} $$ 

One readily sees that $ \phi_N $ is a quasiisomorphism for 
$ N =  \Homder_{\cal A }^{\bullet}( K^{\bullet}, K^{\bullet}_w )$. But these objects
generate $ {\cal D}^b( {\Mof-}  \Endder_{\cal A }^{\bullet}( K^{\bullet} )) $ 
and we can argue as above.

\medskip
\noindent

Let us finally comment on the shifts $ [1] $. It should here be noticed that 
since we are working with right DG-modules, the appropriate definition of 
$ { \cal M }[1] $ for a DG-module $ \cal M $ over a DG-algebra $ \cal A $ 
is the following: 
$$ ( M[1])^i = M^{i+1}, \, \,  d_{ M[1] } = -d_M  \, \, \, { \rm and } 
\, \, \,  m \circ a = ma $$ 
where $  m \circ a $ is the multiplication in $ M[1] $ while $ ma $ is the
multiplication in $ M $. So unlike the left module situation, the $ A $ 
structure on $ M $ is here left unchanged. One now checks that the functors 
commute with the shifts $ [1] $.

\end{pf*}

\medskip
\noindent

\medskip

Now $ A $ has a Koszul grading so we may consider the category
${\mof-}A $ of finitely generated, graded $ A $-modules.
Let us use a lower index to denote 
the graded parts with respect to this $ \mathbb Z $-grading.
This passes to the derived category of graded DG-modules 
which we denote by
$ {\cal D}^b({\mof}- {\cal A} ) $. It may be identified 
with the derived category of ${\mof-}A  $ and  carries a 
twist $ \langle 1 \rangle $ which we arrange the following way:
$$ M \langle 1  \rangle_i = M_{i-1} $$ 

Now the grading on $ K^{\bullet} $ 
induces a grading on 
$ \Endder_{\cal A }^{\bullet}( K^{\bullet} ) $ as well; it satisfies the rule
$$  \Endder_{ A }^{\bullet}( K^{\bullet} )_i := 
    \smallend_{ A }^{\bullet}( K^{\bullet},K^{\bullet}\langle -i \rangle  ) $$

\medskip
With this grading on $ \Endder_{\cal A }^{\bullet}( K^{\bullet} ) $, we may 
consider its graded module category, which we denote
$ {\cal D}^b ( {\mof}- \Endder_{\cal A }^{\bullet}
\left( K^{\bullet}) \right)$.
Let us moreover denote the subcategory of 
this generated, with twists, by the graded summands of the DG-module
$  \Endder_{\cal A }^{\bullet}( K^{\bullet}) $ itself by 
$$ {\cal D}^b ( \langle \Endder_{\cal A }^{\bullet}\left( K^{\bullet}) \rangle  
\right)$$ 
(It would have been more consistent with the earlier notation to write 
$ {\cal D}^b ( {\mof}- \Endder_{\cal A }^{\bullet}
\left( K^{\bullet}) \right)$
for this category;
we however 
prefer the above notation in order to save space). 

\medskip

It is now straightforward to prove the following strengthening of the former
theorem:
\begin{Thm}
The functor 
$ \RHom_{\cal A }(K^{\bullet}, - )$  
establishes an equivalence of the triangulated categories 
$ {\cal D}^b({\mof}- {\cal A} ) $ and
$ {\cal D}^b ( \langle \Endder_{\cal A }^{\bullet}\left( K^{\bullet}) \rangle  
\right)$. It commutes with the twists $ [1]$ and $ \langle 1 \rangle $. 

\end{Thm}

\medskip
\medskip

Our next task will be to study the DG-algebra 
$   \left(  \Extder_{  \cal A }^{\bullet}
( K^{\bullet} ), d= 0  \right)  $. More precisely, we are going to
compare
$  {\cal D}^b \left(   \langle \Extder_{  \cal A }^{\bullet}
( K^{\bullet} ), d= 0 \rangle \right)  $ with 
the standard derived 
category of the category of finitely generated, graded modules over 
\hfill\break
$\Extder_{A}^{\bullet}(k,k)$, i.e.
$ {\cal D}^b \left( {\mof-}\Extder_{ A }^{\bullet}(k,k ) \right) $.
Here, we define the  $ \mathbb Z $-grading on \hfill\break 
$(\Extder_{\cal A}^{\bullet}(K^{\bullet}), d= 0)$ to be 
the negative of the DG-grading. 

\medskip

An object
of $ {\cal D}^b \left({\mof-} \Extder_{ A }^{\bullet}(k,k )\right) $ 
is a bounded complex 
$$ 0 \rightarrow  P^{-n} \rightarrow P^{-n+1} \rightarrow 
\ldots \rightarrow P_0 
 \rightarrow \ldots 
\rightarrow  P^{m-1} \rightarrow P^{m} \rightarrow 0 $$
of projective 
graded $\Extder_{ A }^{\bullet}(k,k ) $-modules 
(so the differential has degree $ 0 $). 
Using this, one constructs a DG-module $\cal  P $ of the DG-algebra 
$ ( \Extder_{  A }^{\bullet}(k,k), d=0 ) $ in the following way:
$$ 0 \rightarrow  P^{-n} \langle -n \rangle \rightarrow 
 P^{-n+1} \langle -n+1 \rangle  \rightarrow \ldots 
\rightarrow P^0 \langle 0 \rangle 
\rightarrow  \cdots P^{m} \langle m \rangle \rightarrow 0. $$ 
In other words, $ {\cal P} $ is the module  
$ {\cal P} = \bigoplus_n P^n $ whose $k$'th DG-piece is
$  {\cal P}^k = \bigoplus_i P^i_{k-i}. $ 
We make $ {\cal P} $  
into a graded DG-module by the rule 
${\cal P}_j = \bigoplus_i P^i_{-j} $. 
This construction defines a functor
$$ F: \, \, 
{\cal D}^b \left( {\mof-}\Extder_{ A }^{\bullet}(k,k )\right) \rightarrow
{\cal D}^b \left( {\modc-} (\Extder_{ A }^{\bullet}(k,k ), d=0) \right) $$
since a graded homotopy between two morphisms $ f,g: P^{\bullet} \rightarrow 
Q^{\bullet} $ of modules will be mapped to a homotopy in the DG-category.

\medskip

I will now argue that actually $ F $ is a functor into the category
\hfill\break 
$ {\cal D}^b \left( \langle \Extder_{   A }^{\bullet}(k,k),d=0 \rangle \right) $. To see
this note first that the $  P_i \langle i \rangle $ all lie in this 
category since they are summands of (shifts of) 
$ ( \Extder_{  A }^{\bullet}(k,k), d=0 ) $. On the other hand $ F $ can be viewed 
as an iterated graded mapping cone construction in the DG-sense, starting 
with the 
morphism 
$ P^{m-1} \langle m \rangle  \rightarrow P^{m} \langle m \rangle $,
which is already a DG-morphism, 
and the claim follows. 
The graded mapping cone
of a graded DG-module morphism 
$ {\cal M} \stackrel{u}{\rightarrow} {\cal N} $ is given by 
$ C(u)= N \oplus M[1] $ 
as DG-module, and reversing the degrees with respect to the 
$ \mathbb Z $-grading.

\medskip

It is now clear that $ F $ commutes with the twists $ [1] $. On the other 
hand $ F $ takes a sequence of 
$ \Extder_{  A }^{\bullet}(k,k) $-modules isomorphic to a standard triangle 
$$ M \stackrel{u}{\rightarrow}  N \rightarrow C(u) 
\rightarrow M[1] $$
to a sequence of graded DG-modules isomorphic to 
$$ F(M) \stackrel{F(u)}{\rightarrow} F(N) \rightarrow F(C(u)) 
\rightarrow F(M[1]). $$

But this is easily seen to be a standard triangle of DG-modules; in other 
words $ F $ is a triangulated functor.
\medskip

One furthermore checks that $ F $ is full and faithful and it is thus an 
equivalence of triangulated categories once we have shown that the generators
$ ( \Extder_{  A }^{\bullet}(k, k_w), d=0 ) $ lie in the image of $ F $. But this is clear.

\medskip

Since 
$$ M_1 \rightarrow M_2 \rightarrow M_3 \rightarrow $$ is a triangle if and
only if 
$$ F(M_1) \rightarrow F(M_2) \rightarrow F(M_3) \rightarrow $$ is a triangle
we deduce that the inverse functor $ G $ is triangulated as well.

\medskip

Our next task will be to analyze the behavior of the twists 
$ \langle 1 \rangle $ with respect to $ F $. Now apart from the change of 
the sign of the differentials, we obtain the relation:
$$ F(M \langle 1 \rangle ) = F(M)[-1]\langle - 1 \rangle. $$
On the other hand, let for a complex $ P^{\bullet} $ of graded 
$ \Extder_{  A }^{\bullet}(k,k) $-modules $ P = \bigoplus P^k_n $ be 
the total module. We may then 
consider the map $ \sigma $ on $ P $ defined as follows:
$$  p \in  P^k_n \mapsto (-1)^n p. $$ 
One now checks that $ \sigma $ defines an isomorphism between the functors 
$ F( - \langle 1 \rangle ) $ and $ F(- )[-1]\langle -1 \rangle. $

\medskip
\medskip

Recall once again that since $ A \, $ is Koszul, the DG-algebras 
$$ \Endder_{\cal A }^{\bullet}\left(K^{\bullet}\right)\mbox{ and }
(\Extder_{ A }^{\bullet}(k,k), d=0)$$ are 
quasiisomorphic so there is an equivalence
$$ \phi: \,  
{\cal D}^b ( {\Modc-}  \Endder_{\cal A }^{\bullet}\left( K^{\bullet} )\right) 
\rightarrow
{\cal D}^b \left( {\Modc-}(\Extder_{ A }^{\bullet}(k,k ),d=0) \right). $$ 
A closer  
look at the above quoted proof reveals that the  quasiisomorphism is 
even a graded one so we obtain an equivalence
$$ \phi: \,  
{\cal D}^b ( \langle  \Endder_{\cal A }^{\bullet}\left( K^{\bullet}) \rangle
\right) 
\rightarrow
{\cal D}^b \left( \langle  \Extder_{ \, A }^{\bullet}(k,k ),d=0  \rangle
\right) $$

\medskip 

Let us gather the results in one theorem.
\begin{Thm} The composition of the above functors
$$ G \circ \phi \circ \RHom_{\cal A }(K^{\bullet}, - ):
{\cal D}^b ( {\mof-} A ) \rightarrow
{\cal D}^b \left( {\mof-}\Extder_{ A }^{\bullet}(k,k )\right) $$ is an 
equivalence of triangulated categories: the Koszul duality functor.
It commutes with $ [1] $ and satisfies the rule
$$ F(M \langle 1 \rangle ) = F(M)[-1]\langle -1 \rangle. $$
\end{Thm}

\section{\bf Translation and Zuckerman functors.}

\noindent
Let us return to the translation functors to and from the wall
$ T_0^{\lambda}:\, \cal{O}\, \rightarrow \cal{O}_{\lambda} $,
$ T^0_{\lambda}:\, \cal{O}\, \rightarrow \cal{O}_{\lambda} $. We saw in 
the first section, that they can be viewed as functors between the 
categories $ \Modc -A $ and $ \Modc -A_{\lambda} $. Now, in [2] it is 
shown that they can be lifted to graded functors 
$$ T_0^{\lambda}:\, \modc -A \rightarrow \modc -A_{\lambda} $$
$$ T^0_{\lambda}:\, \modc -A_{\lambda} \rightarrow \modc -A  $$
that still are adjoint and exact and such that $ T_0^{\lambda} $ takes 
pure objects of weight $ n $ to pure objects of the same weight;
this construction is based on 
Soergel's theory of modules over the coinvariants of the Weyl group. See also 
[1], where graded translation functors are constructed in the setting of 
modular representation theory. 

\medskip
\noindent
The graded translation functor $ T_0^{\lambda} $ induces a homomorphism of 
graded DG-algebras 
$ \Endder_{\cal A }^{\bullet}( K^{\bullet} ) \, \rightarrow 
\Endder_{\cal A_{\lambda} }^{\bullet}(T_0^{\lambda}K^{\bullet} ) $. But 
$ T_0^{\lambda} $ and $ T^0_{\lambda} $ are exact and 
$  K^{\bullet} $ is a projective graded resolution of the pure module $ k $, 
so $ T_0^{\lambda}K^{\bullet} $ is a projective graded resolution of 
$ k_{\lambda} $. Hence
$ T_0^{\lambda}K^{\bullet} $ is graded homotopic to the Koszul complex
$ K^{\bullet}_{\lambda}$ of 
$ A_{\lambda} $ and this implies that 
$ \Endder_{\cal A_{\lambda} }^{\bullet}(T_0^{\lambda}K^{\bullet} )$ and 
$ \Endder_{\cal A_{\lambda} }^{\bullet}(K^{\bullet}_{\lambda} )$ are 
quasiisomorphic graded DG-algebras. 

\medskip

We already saw that $ A $ being a Koszul ring implies that
$$ \Endder_{\cal A }^{\bullet}( K^{\bullet} ) \cong
( \Extder_{\cal A }^{\bullet}( k,k ),d=0 ) $$ 
and by the above and since $A_{\lambda}$ is Koszul as well we have that
$$ \Endder_{\cal A_{\lambda} }^{\bullet}(T_0^{\lambda}K^{\bullet} ) \cong
(\Extder_{\cal A_{\lambda} }^{\bullet}( T_0^{\lambda} k,T_0^{\lambda} k),d=0 ) 
\cong (\Extder_{\cal A_{\lambda} }^{\bullet}(k_{\lambda}, k_{\lambda}),d=0 )$$
We now obtain the following commutative diagram of functors:
$$\hss\begin{array}{ccccc}
{\cal D}^b({\modc-}A) & 
\stackrel{{\scriptscriptstyle - \, \otimes \, K^{\bullet}}}{\leftarrow}
& {\cal D}^b(\langle \Endder_{\cal A}^{\bullet}(K^{\bullet})\rangle) & 
\leftarrow &  {\cal D}^b(\langle \Extder_{\cal A}^{\bullet}(k,k), d=0 \rangle)\\
\downarrow {\scriptscriptstyle T_0^{\lambda}} & & 
\downarrow
& & 
 \downarrow  {\scriptscriptstyle T_0^{\lambda} } \\
{\cal D}^b({\modc-}A_{\lambda}) & 
\stackrel{{\scriptscriptstyle - \,\otimes \,  T_0^{\lambda}
K^{\bullet}}}{\leftarrow}
& {\cal D}^b(\langle \Endder_{\cal A_{\lambda}}^{\bullet}(T_0^{\lambda}
 K^{\bullet} )\rangle) &
\leftarrow &  {\cal D}^b
(\langle \Extder_{\cal A_{\lambda}}^{\bullet}
(T_0^{\lambda}k,T_0^{\lambda}k), d=0 \rangle)
\end{array}
\hss$$
where the middle vertical arrow is the graded tensor product functor
$$ - \otimes_{ \Endder_{\cal A }^{\bullet}
( K^{\bullet} )}\Endder_{\cal A }^{\bullet}( T_0^{\lambda} K^{\bullet}) $$
while the first and third vertical arrows are the graded translation functors. 
The commutativity of the first square is here obvious, whereas the 
quasiisomorphism of $ ( \Extder_{\cal A }^{\bullet}( k,k ),d=0 )) \times 
\Endder_{\cal A_{\lambda} }^{\bullet}(T_0^{\lambda}K^{\bullet} ) $-bimodules 
$$ \Endder_{\cal A }^{\bullet}( T_0^{\lambda} K^{\bullet}) \cong
 \Extder_{\cal A_{\lambda} }^{\bullet}(T_0^{\lambda}k,T_0^{\lambda }k) $$ 
gives the natural transformation that makes the second square commutative.

\medskip

By the same reasoning as in the last section, the two lower arrows define
equivalences of triangulated categories.

\medskip

Now Backelin [2] shows that one can choose the isomorphisms 
$ \Extder_{\cal A }^{\bullet}( k,k ) \simeq A $ and
$ \Extder_{\cal A_{\lambda} }^{\bullet}(T_0^{\lambda} k,T_0^{\lambda} k) \simeq 
A^{\lambda} $ to obtain the following commutative diagram:
$$ \begin{array}{ccc}
   \Extder_{\cal A }^{\bullet}( k,k) & \stackrel{\sim}{\longrightarrow}  & A \\
   {\scriptscriptstyle T^{\lambda}_0}  \downarrow  & & \, \,\, \, \, 
   \downarrow {\scriptscriptstyle \tau}  \\
   \Extder_{\cal A_{\lambda} }^{\bullet}(T_0^{\lambda} k,T_0^{\lambda} k)  
   &  \stackrel{\sim}{\longrightarrow} &  A^{\lambda} 
   \end{array}
$$
The key point is here that both vertical maps are surjections: 
for $ \tau $ this 
is clear, while for $ T^{\lambda}_0 $ an argument involving the Koszul
property of $ A $ is 
required (one can here give a simple alternate argument along the lines 
of the Cline, Parshall, Scott approach to Kashdan-Luzstig theory [8]). 
It then follows that the kernels
of the two vertical maps are the ideals generated by corresponding idempotents
(thus in the degree $ 0 $ part). 

Although the diagram involves non-graded maps, it can be used to give $ \tau $
a grading -- and then of course it is a commutative diagram of graded 
homomorphisms.

\medskip

This diagram, on the other hand, gives rise to the following commutative
diagram of functors:
$$ \begin{array}{ccc}
  {\cal D}^b  (\langle \Extder_{ \cal A }^{\bullet}( k,k),d=0 \rangle ) & 
\stackrel{F}{\leftarrow}    & 
{\cal D}^b({\modc-} A) \\
   {\scriptscriptstyle T^{\lambda}_0}  \downarrow  & & \, \,\, \, \, 
   \downarrow {\scriptscriptstyle \tau}  \\
  {\cal D}^b  ( \langle  \Extder_{ \cal A_{\lambda} }^{\bullet}
(T_0^{\lambda} k,T_0^{\lambda} k),d=0 \rangle ) 
   &  \stackrel{F}{\leftarrow} &  {\cal D}^b ( {\modc-}A^{\lambda})  
   \end{array}
$$

\medskip

We now join the two diagrams of functors to obtain a diagram, in which 
the upper arrows compose to the Koszul duality functor of $ \cal O $ while  
the composition of the lower arrows is isomorphic to the Koszul duality 
functor of $ \cal O_{\lambda} $. Let us formulate this as a Theorem

\begin{Thm}
The translation-- and Zuckerman functors are Kozsul to each other, in other 
words there is a commutative diagram of functors:
$$ \begin{array}{ccc}
  {\cal D}^b ( {\modc-}A ) &  
\stackrel{D}{\leftarrow}    & 
{\cal D}^b({\modc-} A) \\
   {\scriptscriptstyle T^{\lambda}_0}  \downarrow  & & \, \,\, \, \, 
   \downarrow {\scriptscriptstyle \tau}  \\
  {\cal D}^b  {(\modc-}A_{\lambda} ) & 
 \stackrel{D_{\lambda}}{\leftarrow} &  {\cal D}^b ( {\modc-}A^{\lambda})  
   \end{array} $$
where $D$ (resp. $D_{\lambda} $) is the Koszul duality functor as 
described above. 

\end{Thm}

This is the Theorem announced in the introduction of the paper.

\section{\bf The Enright-Shelton Equivalence }

\noindent
We now consider the categorification of the Temperley-Lieb algebra. Let 
thus $ \mathfrak g :=  gl_n $ and $ {\cal O := \cal O} ({ \mathfrak gl}_n )  $.
Let $ \epsilon_k, k \in \{1,2, \ldots, n \} $ be the standard basis of the 
weight lattice and define for $ k \in \{1,2, \ldots, n \}, \lambda_k := 
\epsilon_1 + \epsilon_2 + \ldots \epsilon_k $. We then let 
$ {\cal O}_{k,n-k} $ be the singular
block of $ \cal O $ consisting of modules with central character 
$ \theta ( \lambda_k ) $. So the Verma module 
with highest weight $ \lambda_k - \rho $ lies in $ {\cal O}_{k,n-k} $. Let 
$ {\mathfrak g}_i, 1 \leq i \leq i-1 $ be the subalgebra of $ \mathfrak g $ 
consisting of the matrices whose entries are nonzero only on the 
intersection of the $ i $-th and $ ( i+1 ) $-th rows and columns. We then 
denote by  $ {\cal O}_{k,n-k}^i  $ the parabolic subcategory of 
$ {\cal O}_{k,n-k} $ whose modules are the locally ${\mathfrak g}_i $-finite
ones in $ {\cal O}_{k,n-k} $.

\medskip

We shall also consider the following dual picture: let $ { \mathfrak p }_k $ 
be the parabolic subalgebra of $ { \mathfrak g } $, whose Levi part is 
$ { \mathfrak g }_k \oplus  { \mathfrak g }_{ n-k } $ and such that 
$ { \mathfrak n }_+ \subseteq  { \mathfrak p }_k $ and let 
$ {\cal O}^{k,n-k} $ be the the full subcategory of $ \cal O $ consisting of 
the locally $ { \mathfrak p }_k $-finite modules. Choose an integral
dominant regular weight $ \mu $ and integral dominant subregular weights 
$ \mu_i $ on the $ i $-th wall, $ i = 1, 2, \ldots , n-1 $ (so the
 coordinates 
of $ \mu_i $ in the $ \epsilon_i $ basis are strictly decreasing, except for 
the $ i $-th and the $ ( i+1 )$-th that are equal). Let finally 
$ {\cal O}^{k,n-k}_i  $ be the subcategory of $ {\cal O}^{k,n-k} $ with 
central character $ \theta ( \mu_i ) $. 

\medskip
All of this is the setup of [6]

\medskip

Let $ R^i_{k,n-k} ( $resp.$  R_i^{k,n-k}  $) be the endomorphism ring of 
the minimal projective generator of $ {\cal O}_{k,n-k}^i $  
(resp. $ {\cal O}^{k,n-k}_i $). The following theorem is a direct 
consequence of Backelin's work [2]:

\begin{Thm}
$ ( R_{k,n-k}^i)^! = {R}^{k,n-k}_i $
\end{Thm}
\noindent
\begin{pf*}{Proof} 
The main result of [2] is that $$ ( R^{\psi}_{\phi} )^! =  
R^{\phi}_{-w_0 \psi} $$ with the notation as in [5]. 

Now one 
should first observe 
that $ c \, Id \in  { \mathfrak gl}_n  $ acts on 
$ { \cal O }_\lambda $ through multiplication with 
$ c  \sum  \lambda_i $; hence $ {\cal O}_\lambda $ is equivalent to 
the category $ {\cal O}_{\overline{\lambda}} $ of 
$ { \mathfrak sl}_n  $-modules, 
where $\overline{\lambda} $ denotes the image of $ \lambda $ under 
the projection
of the weight lattice with kernel $ {\mathbb Z} \sum \epsilon_i $. We can 
thus restrict ourselves to the semisimple situation and may indeed use the
results quoted. 

Now $ \mu_i $ and $ { \mathfrak sl}_i \subseteq 
{ \mathfrak sl}_n $ are both given by the simple root $ \alpha_i = 
\epsilon_i -\epsilon_{i+1} $ and also $ { \mathfrak p}_{k,n-k} $ and 
$ \lambda_k $ are given by the same simple roots ($ = \Delta \setminus
\{ \alpha_i \} $). Finally, we have in type A that 
$ w_0 \lambda = - \lambda $ and the theorem follows.
\end{pf*}

As a corollary we obtain a simple proof of the following equivalence of 
categories first proven by Enright and Shelton.

\begin{Corollary}
$ {\cal O}^{k,n-k}_1   \cong {\cal O}^{k-1,n-k-1} $ 
\end{Corollary}

\begin{pf*}{Proof} 
There is first of all a standard equivalence of categories: 
$$ {\cal O}_{k,n-k}^1     \cong {\cal O}_{k-1,n-k-1} $$
The functor $ \nu_n $ from left to right takes the sum of all 
weight spaces of weight $ \epsilon_1 + x_3 \epsilon_3 + \ldots 
x_n \epsilon_n - \rho_n $ with $ x_i \in \mathbb Z $. This is a 
$ { \mathfrak gl}_{n-2} $-module and $ \nu_n $ is then the tensor product of
it with the module defined by 
$ \epsilon_1 +  \epsilon_2 + \ldots \epsilon_{n-2} $. The inverse functor 
comes from an induction procedure, see [6] for the details.
Now this equivalence gives us a ring isomorphism 
$$  (R_{k,n-k}^1)^! \cong  ({R}_{k-1,n-k-1})^! $$ which 
combined with the theorem 
yields a ring isomorphism 
$$ \xi:  R^{k,n-k}_1 \cong {R}^{k-1,n-k-1} $$ 
But then the module categories of $ R^{k,n-k}_1 $ and $ {R}^{k-1,n-k-1} $ 
are equivalent, by the restriction and extension of scalars along $ \xi $.
The corollary is proved.
\end{pf*}

\medskip

\noindent
This might be useful in proving the full conjectures of [6].

\end{document}